\documentclass{amsart}
\usepackage{amssymb,latexsym,amsmath}

\theoremstyle{plain}
\newtheorem{theorem}{Theorem}[section]
\newtheorem{corollary}[theorem]{Corollary}
\newtheorem{lemma}[theorem]{Lemma}
\newtheorem{proposition}[theorem]{Proposition}
\newtheorem{definition}[theorem]{Definition}

\newtheorem{remark}[theorem]{Remark}

\numberwithin{equation}{section}

\font\cmcsc=cmcsc10 at 8pt

\begin{document}

\title[bare term structure]
        {One-Factor Term Structure without Forward Rates}
\author{Victor Goodman}
\address{Mathematics Department\\
         Indiana University\\
         Bloomington, IN 47405}
\email{goodmanv@indiana.edu}
\author{Kyounghee Kim}
\address{Mathematics Department\\
         Indiana University\\
         Bloomington, IN 47405}
\email{kyoukim@indiana.edu}

\begin{abstract}
We construct a no-arbitrage model of bond prices where the long bond is used as a numeraire.  We develop bond prices and their dynamics without developing any model for the spot rate or forward rates.  The model is arbitrage free and all nominal interest rates remain positive in the model.  We give examples where our model does not have a spot rate; other examples include both spot and forward rates.
\hfill\break
\phantom{as}\hfill\break
\noindent {\cmcsc Keywords}:\ term structure, no-arbitrage, zero coupon bonds, stochastic volatility \hfill\break
\noindent {\cmcsc Subject Classification}:\   Primary 91B28, 60H30\ \      Secondary 60J65
\end{abstract}

\maketitle

\vspace{4ex}

\section{introduction}

Models of term structure fall into two general categories: equilibrium models and no-arbitrage models.  All the current models in each of these categories include spot rates and forward interest rates as part of the model structure (see  [Mi] and [MR]).  We introduce a no-arbitrage model for term structure that does not incorporate such rates.   Although our model fits any reasonable data for the initial term structure exactly, some lattitude remains for a modeler to interpolate initial data so that interest rates may appear within the model.  This option is a new feature for model calibration which might be exploited; for instance, implied interest rates may be developed from interest rate data to better fit the additional data.  An example of this is given at the end of Section 6.

We define a one-factor term structure model using a standard Brownian motion process which we denote by   $B_t$. Our basic stochastic processes will be two well-known functions of Brownian motion:

 $M_t$ will denote the simple exponential martingale
$$M_t = \exp(B_t - \frac{t}{2})$$
and $A_t$ will denote its time integral
$$A_t = \int_0^t M_sds.$$
Our bond prices will be defined in terms of the processes $M_t$ and $A_t$, and the dynamics of these asset prices will be determined by the dynamics of  these two processes.  A scaling property of Brownian motion allows us to generalize our model in a straightforward way to include a volatility parameter (See Section 5).  

Our model places a special emphasis on a finite time horizon.  Following the customary modeling approach,  we calibrate our model at time $t=0$  with bond price data, the market price of discount bonds of various maturities.   We  use a standard notation (see \cite{HJM}, for example) for these observed market prices,  $P(0,T)$, the parameter $T$ being the maturity date.  At maturity each bond is worth one monetary unit.  The horizon date $T_h$ is the maturity date of the long bond.

\vspace{3ex}

\begin{definition}\label{horiz}   An horizon date $T_h$ is any fixed positive number for which there is an observable `calibration value' for $P(0,T_h)$.  

The market quote is called the \emph{long bond price} at time zero and the model price for the \emph{long bond} at any time $t\leq T_h$  is denoted by $P(t,T_h)$.

\end{definition}

\vspace{3ex}

The stochastic process for the long bond is a crucial part of our model.  This maturity sets the longest allowable maturity within the model as well as the time horizon for the evolution of the bond prices. Moreover,  since the value at time $t$  of our numeraire is  $P(t,T_h)$, it is essential to define this process.

However, one probabilistic feature of this model allows us to immediately describe other maturity prices once the long bond process is set.  This description involves a minimal amount of calibration data (whereas the long bond price involves a maximal amount of data) Therefore,  we now describe, indirectly,  the other maturity prices with the following `out-of-order' definition:

\vspace{3ex}

\begin{definition} \label{dis}  For any maturity date  $T <T_h$  for which there is an observable `calibration value' for $P(0,T)$, the stochastic process for the discounted T-bond price is given by
\vspace{1ex}

\begin{equation}\label{discount}
\frac{P(t,T)}{P(t,T_h)} := 
\Big(\frac{P(0,T)}{P(0,T_h)}\Big)^{\big[\frac{2M_t}{2+A_t\log [P(0,T) /P(0,T_h)]}\big]}
\end{equation}

\end{definition}
\vspace{3ex}

\begin{remark}  The stochastic process defined in equation (\ref{discount}) is meaningful precisely when the initial data satisfies $P(0,T) \geq P(0,T_h)$.  Otherwise, a negative logarithm in the definition allows a collapse of the bond price to zero.  This is a restriction on our model to \emph{reasonable calibration data}. 

The initial discount bond prices must be monotone in the maturity date.  This restriction, of course, does allow inversions of the yield curve, which is sometimes seen in real data.  The restriction rules out only \emph{negative} nominal forward rates within the calibration data.

\end{remark}

Equation (\ref{discount}) will be used to define a T-maturity bond price once the long bond price process is chosen.  Notice that the only data dependence within this stochastic process is the (constant)  ratio of the initial bond prices.  From  Definition \ref{dis} one sees that this process has continuous, positive sample paths and the initial value of the process matches the ratio of the market bond prices at time $0$.  At this point we are able to illustrate our choice of the dynamics with the following proposition:

\vspace{3ex}

\begin{proposition}\label{arbitrage} The process given by equation (\ref{discount}) is a positive local martingale iff $P(0,T) \geq P(0,T_h)$.

\end{proposition}
\begin{proof}  Let $F = \log (P(0,T) / P(0,T_h))$ and assume that $F > 0$.  Equation (\ref{discount}) can be written as

$$ \frac{P(t,T)}{P(t,T_h)} = \exp \Big(\frac{2M_t}{2F^{-1}+A_t}\Big) $$

\noindent We let

$$Y_t  = \frac{2M_t}{2F^{-1}+A_t}$$
and a basic calculation shows that 

\begin{equation}\label{sde}
dY = YdB - \frac{1}{2}Y^2dt
\end{equation}
\noindent That is,

$$Y_t = F + \int_0^t \{Y_sdB_s -  \frac{1}{2}Y_s^2ds\}$$
\vspace{2ex}

\noindent This shows that $Y_t$ is the exponent of a \emph{stochastic exponential} (see [WH] or [KS] for a discussion of stochastic exponentials and their relation to the Girsanov theorem).  Hence, $\exp(Y_t)$ is a positive local martingale.  In the case that $F < 0$, $Y_t$ has an initial negative value and explodes to $-\infty$.  In this case, the process in equation (\ref{discount}) is not a positive local martingale.
\end{proof}

\section{ calibrating the long bond}

\vspace{3ex}

Calibration data typically includes a finite sequence of bond prices corresponding to maturity dates $T_1 < T_2 < ... < T_h$.  In fitting a model to this sequence, one frequently interpolates data and effectively determines an artificial bond price $P(0,T) $ for all $T\leq T_h$.  For instance, one plots yields to maturity for the data at hand, and fits a yield curve through finitely many points.  Then additional, artificial bond prices are determined from the completed yield curve.

Some method of interpolating the observed bond prices must be carried out in order to define the stochastic process for the long bond price, $P(t,T_h)$.  Each method of interpolation gives arise to a different dynamics for the long bond.  This in turn affects the dynamics of each shorter maturity (although it has no effect on the stochastic behavior of the discounted bonds in the model, as Proposition \ref{arbitrage} shows).  

Possible benefits of an interpolation scheme will be discussed in Section \ref{spotrate}.  For now, we assume that some method has produced a monotone decreasing function $P(0,T)$ which passes through each observed bond price given by data. The simplest choice is the piece-wise linear curve through these data points.

\vspace{3ex}

\begin{definition} \label{long}  Given an horizon date  $T_h >0$  and a continuous, monotone decreasing function $T \to P(0,T)$ of `calibration values' $, 0 \leq T \leq T_h$ such that $P(0,0)=1$, the stochastic process for the long bond price is given by
\vspace{1ex}

\begin{equation}\label{longb}
P(t,T_h) := 
\Big(\frac{P(0,t)}{P(0,T_h)}\Big)^{\big[\frac{-2M_t}{2+A_t\log [P(0,t) /P(0,T_h)]}\big]}
\end{equation}

\end{definition}
\vspace{3ex}

\begin{remark}\label{longre}  The stochastic process defined in equation (\ref{longb}) has positive, continuous sample paths on the time interval $t \leq T_h$.  Notice that the exponential base is a variable whose continuity properties depend explicitly on the continuity of initial bond prices as a function of their \emph{maturity dates}.  Moreover, it is clear from equation (\ref{longb}) that
$$\lim_{t \downarrow 0}P(t,T_h) = P(0,T_h)$$ 

\vspace{2ex}

Although one might choose the initial curve of maturities to be smooth, the process in equation (\ref{longb}) does not inherit this smoothness because of the stochastic exponent.  Nevertheless, as $t\to T_h$ the bond price converges to one, as one would expect in a term structure model.  

\vspace{2ex}

Moreover, since the variable base is greater than one and the exponent is negative, the long bond price is always less than one before maturity.  This is a feature that many models currently used in term structure modeling do not have (see [ MR, chapters 12 and 15] ).  In our model, the yield to maturity is always positive.

\vspace{2ex}

Definition \ref{long} is motivated by the discounted price expression in equation (\ref{discount}).  If one formally sets the maturity parameter $T$ equal to the time parameter $t$ in  (\ref{discount})  and imposes the value $P(t,t)=1$, one obtains the reciprocal of the price given in Definition \ref{long}.

\end{remark}

\vspace{3ex}

  We use the formula appearing in Definition \ref{dis} and the long bond formula (\ref{longb}) to define all the basic asset prices in the model:

\vspace{3ex}

\begin{definition}\label{asset} The stochastic process for each T-maturity bond price, $T\leq T_h$ is

\begin{equation}\label{assetb}  
P(t,T)= P(t,T_h)\Big(\frac{P(0,T)}{P(0,T_h)}\Big)^{\big[\frac{2M_t}{2+A_t\log [P(0,T) /P(0,T_h)]}\big]}
\end{equation}

\noindent The time parameter $t$ varies over the interval $[0, T]$.
\end{definition}

\vspace{3ex}

\begin{proposition}  For each maturity date $T\leq T_h$, the T-maturity bond price has positive, continuous sample paths.  Moreover, each price process exhibits continuous, limiting price behavior  at its calibration time and at its maturity date.  That is,

$$  \lim_{t \downarrow 0}P(t,T) = P(0,T)  \mbox{\hskip.2 in and \hskip.2in  } \lim_{t \to T}P(t,T) = 1 $$

\end{proposition}

\vspace{3ex}

\begin{proof} Arguments for the continuity behavior in the case $T=T_h$ were given in Remark \ref{longre}.  For $T < T_h$, the product in equation (\ref{assetb}) is clearly positive and continuous.  First, we have 
$$  \lim_{t \downarrow 0}P(t,T) =  \lim_{t \downarrow 0}P(t,T_h)  \frac{P(0,T)}{P(0,T_h)}=P(0,T)
$$
Secondly we observe from Definition \ref{long} that
$$ P(T,T_h) = \lim_{t \to T}P(t,T_h)$$
is the expression for the reciprocal of the second factor in equation (\ref{assetb}) evaluated at $t=T$.  This shows that
$$ \lim_{t \to T}P(t,T) = 1 $$

\end{proof}

\begin{remark}  The proposition shows that each asset price evolves continuously from its initial value.  In particular, if some  $T_i$ corresponds to an observed price at the time of calibration, then $P(t,T_i)$ evolves from this price \emph{whatever interpolation scheme is employed to fill out the initial term structure} (subject to the conditions of Definition \ref{long}).  

Of course, equation (\ref{assetb}) defines artificial bond processes for each interpolated value $P(0,T)$.  Remarkably, the time evolutions of these processes do not influence the processes $P(t,T_i)$.  This unusual feature of the model emphasizes the fact that instantaneous forward rates are not an essential feature of the model dynamics.
\end{remark}

\vspace{3ex}

In order to use techniques of stochastic analysis in developing our model, we establish a basic property of the dynamics given implicitly by equation (\ref{assetb}).

\begin{definition}\label{curve}  If a calibration scheme satisfies the conditions of Definition \ref{long}, the initial price curve $T\to P(0,T)$ defines the \emph{nominal forward rate curve} $F(t)$ where

$$ F(t) = \log (P(0,t) / P(0, T_h))$$

\end{definition}

\begin{remark}  $F(t)$ is a continuous monotone decreasing function on $[0,T_h]$ and $F(T_h)=0$.  Consequently, $F(t)$ is a nonnegative function with finite variation.

\end{remark}

\vspace{3ex}

\begin{proposition}\label{eee}

\begin{equation}\label{eeel}
P(t,T_h) = \exp (-\frac{2M_t }{2 F(t)^{-1} + A_t})
\end{equation}

\noindent  and 

\begin{equation}\label{eeeb}
P(t,T) = \exp (\frac{2M_t }{2 F(T)^{-1} + A_t} - \frac{2M_t }{2 F(t)^{-1} + A_t})
\end{equation}

\end{proposition}

\vspace{3ex}

\begin{proof}  We may write equation (\ref{longb}) as
$$\log P(t,T_h) = -2F(t)M_t / (2+ A_t F(t))$$
$$=  -\frac{2M_t }{2 F(t)^{-1} + A_t}$$
and this establishes equation (\ref{eeel}).  From equation (\ref{assetb}) we have
$$\log P(t,T) = \log P(t,T_h) + 2F(T)M_t / (2+ A_t F(T))$$
and equation (\ref{eeeb}) is an immediate consequence.

\end{proof}

\vspace{3ex}

\begin{theorem}\label{positive}  Each asset price process in the model defined by equations (\ref{longb}) and (\ref{assetb}) is a semimartingale. 
 
\vspace{1ex}

\noindent Each discounted asset price
$$\frac{P(t,T)}{P(t, T_h)}$$
is a positive local martingale.

\vspace{1ex}

\noindent If $F(t)$ is a continuous, strictly decreasing \emph{nominal forward rate function} with $F(0) = -\log (P(0,T_h))$ then all nominal forward rates within the model are positive.

%
%

\end{theorem}

\vspace{3ex}

\begin{proof}  Since the function $F(t)$ has finite variation, the processes $F(t)M_t $ and $F(t)A_t$ are semimartingales (see [P, Chapter 2, Theorem 33]).  It also follows from Theorem 33 that any smooth function of these two semimartingales is also a semimartingale. 

Proposition \ref{eee} states both the long bond price and the T-maturity bond price are exponential functions of positive rational combinations of the processes.  Hence, each asset price is a semimartingale. 

The local martingale property of the discounted assets was established in Proposition \ref{dis}.

If $T_1 < T_2$ the \emph{nominal forward interest rate} at $t \leq T_1$ is 
$$\frac{1}{T_2 - T_1}  \log\big(\frac{P(t,T_1)}{P(t, T_2)}\big )$$
Using Proposition \ref{eee} we obtain
$$\log P(t,T_1)  - \log P(t, T_2) = \frac{2M_t }{2 F(T_1)^{-1} + A_t} - \frac{2M_t }{2 F(T_2)^{-1} + A_t}$$
If $F(T_1) > F(T_2)$ this quantity is always strictly positive.  That is, every nominal interest rate within the model is positive.

\begin{remark}  The positivity of nominal rates also follows directly from equation (\ref{discount}).  One may use the equation twice and eliminate the extra term for $\log P(t,T_h)$ to obtain the same formula for the nominal rate.

\end{remark}

\end{proof}

\section{trading strategies and arbitrage}

\vspace{3ex}

We consider tame, self-financing trading strategies, similar to those introduced in [DH, Lemma 2].  

\vspace{3ex}

\begin{definition} \label{tame} If $\ T_1 < T_2 \dots < T_{n-1} < T_h$ are maturity dates, a \emph{trading strategy}
 $\theta$ is a vector
$$\theta = ( \theta_1 , \theta_2, \dots , \theta_n)$$
where each function $\theta_i$ is a bounded, progressively measurable function of $\ t \leq T_h$ and also satisfies $\theta_i(t) = 0$ for $t>T_i$.  
$\theta_i$ represents a position in the $T_i$-maturity bond; $\theta_n$ is the position in the long bond.  Let
\begin{equation*}
T_\theta = \inf \{ t \le T_h\ :\ \theta (t) = 0 \}
\end{equation*}
A trading strategy must satisfy $\theta \equiv 0$ on the (possibly random) interval $(T_\theta, T_h]$.

A strategy $\theta$ is said to be \emph{tame} if  the  gains process
$$ G(t) := \sum_{i=1}^n \int_0^t \theta_i(s)dP(s,T_i)$$
is uniformly bounded below on the interval $[0, T_\theta )$ when \emph{divided by the long bond price $P(t,T_h)$}.  That is,  $G(t \wedge T_\theta )/P(t\wedge T_\theta ,T_h)$ is bounded below by a non-random constant on the interval $[0,T_h]$.

\end{definition}

\vspace{3ex}
\begin{remark}  Since a trading strategy is progressively measurable, the time $T_\theta$ at which all positions are closed out is a stopping time.  We restrict $\theta$ so that $\theta = 0$ for $t > T_\theta$ because the model allows only bond trading. There is no money market in this model.

\end{remark}

\vspace{3ex}

\noindent {\bf Example.}   Let $T_1 < T_2 < T_h$ be maturity dates.  Suppose that $\sigma_i, i = 1,2$, are stopping times with the property that $\sigma_i \le T_i$ almost surely and $\sigma_1 \le \sigma_2$.  We define another stopping time that controls the size of each discounted bond price, $P(t, T_i)/ P(t,T_h)$.  Fix $C > 0$ and let
\begin{equation}\label{stopping}
\tau = \inf \{t\le T_h : \frac{2M_t }{2 F(T_1)^{-1} + A_t}\ge C  \}
\end{equation}
The process in the definition of $\tau$ equals  $\log (P(t, T_1)/ P(t,T_h))$ for $t \le T_1$ but  this process is defined for all $t \le T_h$  and it dominates all other logarithms of discounted prices with later maturity dates.  Therefore, the other stopped process
$$P(t\wedge \tau, T_2)/ P(t\wedge \tau,T_h)  $$
is also bounded.  

Let $a$ and $b$ be any constants and define a strategy where one holds $a$ units of the $T_1$-maturity bond and $b$ units of the $T_2$-maturity bond until the random time $\sigma_1$.  At this time one closes out the position in $T_1$-bonds and invests the proceeds in the $T_2$-bond.  The latter position is closed out at the random time $\sigma_2$. Then
\begin{equation*}
\begin{aligned}
\theta_1(t) = &a 1_{\{t <\sigma_1\}}\\
\mbox{and\hskip.5in}&\\
\theta_2(t) =& b 1_{\{t <\sigma_1\}}
 + (b +  aP(\sigma_1, T_1)/ P(\sigma_1,T_2) )1_{\{\sigma_1 \le t < \sigma_2\}}\\
\end{aligned}
\end{equation*}
The gains process is
\begin{equation}\label{tameex}
\begin{aligned}
G(t) = \int_0^{\sigma_1 \wedge t}a dP(s,T_1)& + \int_0^{\sigma_2 \wedge t}b dP(s,T_2)\\
+\int_{\sigma_1 \wedge t}^{\sigma_2 \wedge t}&a\frac{P(\sigma_1, T_1)}{P(\sigma_1,T_2) }dP(s,T_2)\\
=aP(\sigma_1 \wedge t,T_1)+ bP(\sigma_2 \wedge t,T_2)
&+a\frac{P(\sigma_1, T_1)}{P(\sigma_1,T_2) }[P(\sigma_2 \wedge t,T_2)  - P(\sigma_1 \wedge t,T_2)]\\
\end{aligned}
\end{equation}

However, if either $a$ or $b$ were negative, so that we have a short position in one bond, then the \emph{discounted} gains process would not be bounded below.  Neither of the discounted bond prices $P(s, T_i)/ P(s,T_h) $ is bounded so a constant short position, even up to a nonrandom time,  is not a tame strategy.

On the other hand, the strategy 
$$\theta (t) 1_{\{t \le \tau\}}$$
is tame.  The indicator function modifies the original strategy so that when the time $\tau$ is encountered, all positions are closed out.  One can see from the formula for the gains process in equation (\ref{tameex}) that trading stopped at $\tau$ effectively replaces each time $\sigma_i$ by $\sigma_i\wedge \tau$.  All discounted quantities in the stopped version of equation (\ref{tameex}) are bounded.  Therefore, the modified trading strategy is tame.
\vspace{3ex}

\begin{definition}
A trading strategy $\theta$ is \emph{self-financing} iff the gains process $G(t)$ in Definition \ref{tame} satisfies
$$\theta_t \cdot (P(t,T_1),\dots , P(t,T_h)) = \theta_0\cdot (P(0,T_1),\dots , P(0,T_h)) +G(t)$$
for $t \leq T_\theta$.

\end{definition}

We will show that 
\begin{equation}\label{dismart}
\theta_t \cdot (P(t,T_1),\dots , P(t,T_h)) / P(t,T_h) = 
\end{equation}
$$\theta_0\cdot (P(0,T_1),\dots , P(0,T_h))/P(0,T_h)+\sum_{i=1}^n \int_0^t \theta_i d\big(\frac{P(s,T_i)}{P(s,T_h)}\big)
$$
Equation (\ref{dismart}) is the well-known identity that relates self-financing portfolio values to stochastic integrals of the positions against the discounted bond prices.

Consider the stopped version of equation (\ref{dismart}) where we stop trading at the time $\tau$ defined in equation (\ref{stopping}).  The stochastic integral up to the first time $t $ or $\tau$ is effectively a stochastic integral up to time $t$ with respect to the martingale processes
$$\frac{P(s\wedge \tau,T_i)}{P(s\wedge \tau,T_h)}$$
Each stopped process is a bounded local martingale according to Proposition \ref{arbitrage} and hence is a martingale. Theorem 12.1 of \cite{Do} shows that a discounted portfolio value equals an integral with respect to the discounted asset prices under these hypotheses.  Therefore, the identity we seek holds on the event 
$\{\tau \ge t \}$.  But $\tau$ was defined with a parameter $C$.  As we let $C \to\infty$, the event for which the identity holds expands to an event with probability one.  This establishes equation (\ref{dismart}) for any self-financing trading strategy. 

\begin{lemma}\label{supermg} If $\theta$ is a tame, self-financing trading strategy, then the stopped, discounted portfolio value
\begin{equation*}
\theta_{t\wedge T_\theta} \cdot (P(t\wedge T_\theta,T_1),\dots , P(t\wedge T_\theta,T_h)) / P(t\wedge T_\theta,T_h)
\end{equation*}
is a supermartingale on the interval $[0,T_h]$. 

\end{lemma}

\begin{proof}   Theorem  \ref{positive} states that each process $P(t,T_i) / P(t, T_h)$ is a local martingale.  Since each position process $\theta_i(t)$ is adapted and bounded, each stochastic integral
$$ \int_0^t \theta_i (s)d\big(\frac{P(s,T_i)}{P(t,T_h)}\big)$$
is  a local martingale.  By assumption, the stopped, discounted gains process is bounded below and  it is also a local martingale.  It follows from Fatou's lemma that the discounted gains process, stopped at $T_\theta$, (and therefore the discounted portfolio value as well) is a supermartingale.

\end{proof}

\vspace{3ex}

\begin{lemma}\label{extended} Let 
\begin{equation}\label{hold}
\Pi_{t\wedge T_{cl}}  / P(t\wedge T_{cl},T_h)
\end{equation}
denote a stopped, discounted portfolio value process, where all bond positions are closed out at some stopping time $T_{cl} \le T_h$.  Then the process also equals the discounted value for bond trading which continues to the horizon date.

If trading is extended by investing the monetary value $\Pi_{cl}$ in $T_h$-maturity bonds at time $T_{cl}$ and holding the bonds, then the process in (\ref{hold}) is the discounted value throughout the interval $[0,T_h]$. 

\end{lemma}

\vspace{3ex}

\begin{proof}
At time $T_{cl}$ a monetary  amount $\Pi_{T_{cl}}$ may be exchanged for 
$$ \Pi_{T_{cl}}/ P(T_{cl},T_h) $$
units of the long bond.  Thereafter, the portfolio value equals
$$ \Pi_{T_{cl}}\frac{P(t,T_h)}{ P(T_{cl},T_h)} $$
Hence, the discounted portfolio value is also given by the process in equation (\ref{hold}) for $T_{cl} < t \le T_h$.

\end{proof}

\vspace{3ex}

\begin{theorem}\label{noarbitrage}  The collection of all tame, self-financing trading strategies for the asset prices given by Definitions \ref{long} and \ref{asset} forms a system without arbitrage strategies. 

\end{theorem}

\vspace{3ex}

\begin{proof}   Lemma \ref{supermg} implies that each tame, self-financing strategy produces a portfolio value so that its discounted version, when parametrized by $t\wedge T_\theta$  is a supermartingale on the interval $[0,T_h]$.   In addition, Lemma \ref{extended} implies that  the process is a discounted portfolio value on this time interval.  We apply [Do, Theorem 12.2] which states that if all discounted portfolio values are supermartingales there can be no (tame, extended) arbitrage strategy on the interval $[0,T_h]$.

\end{proof}
\section{examples of interest rate derivatives}

\vspace{3ex}

\noindent {\bf Example 1.} A Forward Contract on the Long Bond
\vspace{1ex}

The obvious trading strategy to replicate a forward contract on a long bond allows us to illustrate a method  for computing prices using local martingales.   We  denote the forward contract price for delivery at time $T < T_h$ by $f(t)$ . The exercise price required to receive a long bond is denoted by $\kappa$. 

If the discounted gains process in equation (\ref{dismart}) were a martingale, then we would identify replicating  portfolio values by computing the right hand expression in the `identity'
\begin{equation}\label{longcontract}
 \frac{f(t)}{P(t,T_h)} = E[ (P(T,T_h) - \kappa)/ P(T,T_h)\ |\  \mathcal{F}_{t} \ ]
\end{equation}
In our case the martingale assumption is incorrect, so we introduce stopping times.  If one ignores the stopping time issue and computes a price using equation (\ref{longcontract}) directly, one gets an incorrect answer (see appendix \ref{pitfall}).

Proposition \ref{eee} implies that
$$Y_t = \log [P(t,T) / P(t,T_h)] = \frac{2M_t }{2 F(T)^{-1} + A_t} $$
and for each $n=1,2, \dots$ we define the random time   $\tau_n$ by
\begin{equation}\label{time}
\tau_n = \inf \{ t \wedge T_h\ ; Y_t \geq n\}
\end{equation}
It follows from an argument in the proof of Proposition \ref{arbitrage} that
$$ \exp ( Y_{t\wedge \tau_n})$$
is a martingale.  We then know that
$$\frac{P(t\wedge \tau_n, T)}{P(t\wedge \tau_n,T_h)}$$
is a martingale.  This suggests another contract that has a payoff similar to the forward, but may terminate early with a different terminal value.  That is, we consider the conditional expectation
\begin{equation}\label{longstopped}
E[ (P(T\wedge \tau_n,T_h) - \kappa P(T\wedge \tau_n,T))/ P(T\wedge \tau_n,T_h)\ |\  \mathcal{F}_{t} \  ]
\end{equation}
$$ = 1 - \kappa E[  \frac{P(T\wedge \tau_n,T)}{P(T\wedge \tau_n,T_h)} \ |\  \mathcal{F}_{t} \ ] $$
We think that the payoff is the long bond minus $\kappa$ shares of the T-bond and the payoff occurs at the smaller of the times $T$ and $\tau_n$. This is useful because we can compute the conditional expected value in equation (\ref{longstopped}).  It is
$$ = 1 - \kappa   \frac{P(t\wedge \tau_n,T)}{P(t\wedge \tau_n,T_h)}  $$
Moreover, there is a tame trading strategy that produces the payoff value
$$ P(T\wedge \tau_n,T_h) - \kappa   P(T\wedge \tau_n,T) $$
It is $\theta(t) = (-\kappa \ , 1)$, for  $ t\leq \tau_n \wedge T$.
Now the payoff value converges to
$$ P(T,T_h) - \kappa  $$
with probability one as $n\to\infty$.  And, the earlier price converges to
$$ P(t,T_h) - \kappa   P(t,T) $$
These limits are portfolio values because the trading strategy `discovered' by computing in the stopped situation also converges to $\theta = (-\kappa \ , 1)\ ,  t\leq T$ ; equation (\ref{dismart})  establishes the link between initial and terminal investment values via a trading strategy, even in the limit as $n\to\infty$.

\vspace{3ex}

\vspace{3ex}

\noindent {\bf Example 2.} A Forward Contract on $T'$-maturity Bond
\vspace{1ex}

In this example the delivery date is $T$ and $T < T'$.  The exercise price required to receive the bond is denoted by $\kappa$. 
To pass to martingale computations we use the same sequence of  stopping times as defined by equation (\ref{time}) in Example 1.

As before,%
$$ \exp ( Y_{t\wedge \tau_n})=\frac{P(t\wedge \tau_n, T)}{P(t\wedge \tau_n,T_h)}$$
is a martingale.  We consider a contingent claim whose payoff is one share of the $T'$-bond minus $\kappa$ shares of the $T$-bond at the smaller of the times $T$ and $\tau_n$.  Equation (\ref{dismart}) implies that a discounted portfolio, stopped at $\tau_n$, is a martingale.  Therefore, we consider
\begin{equation}\label{shortcon}
\begin{aligned}
 E[ (P(T & \wedge \tau_n,T') - \kappa P(T\wedge \tau_n,T))/ P(T\wedge \tau_n,T_h)\ |\  \mathcal{F}_{t} \  ]\\
&=  E[ (P(T\wedge \tau_n,T') / P(T\wedge \tau_n,T_h)\ |\  \mathcal{F}_{t} \  ]\\
&-   \kappa E[  P(T\wedge \tau_n,T)/ P(T\wedge \tau_n,T_h)\ |\  \mathcal{F}_{t} \  ]
\end{aligned}
\end{equation}
Each expected value is the conditional value of a martingale.  The process with maturity date $T$ is a martingale because the stopping time stops the local martingale so that it is a bounded process.  Moreover, this process  dominates the other  local martingale point-wise because the process $Y_t$ is decreasing as a function of the maturity date.  Therefore, the process with the $T'$ maturity date is also a martingale.   The conditional expectations in equation (\ref{shortcon}) equal

$$ = \frac{P(t\wedge \tau_n,T')}{P(t\wedge \tau_n,T_h)}  - \kappa   \frac{P(t\wedge \tau_n,T)}{P(t\wedge \tau_n,T_h)}  $$
so that an initial portfolio value for this contingent claim is
$$  P(t\wedge \tau_n,T')  - \kappa   P(t\wedge \tau_n,T) $$
We see that the trading strategy is $\theta(t) = ( -\kappa, 1, 0)$ for $ t \leq T\wedge \tau_n$. That is, some $T$-maturity bonds are sold short, a $T'$-bond is held, and no long bonds are involved.    As $n\to \infty$ the payoff, the initial portfolio value, and the trading strategy converge in a trivial sense.  Equation (\ref{dismart}) holds on the event $\{Y_t < \infty\}$ which has probability one, therefore we see that the contract price equals
$$  P(t,T')  - \kappa   P(t,T) $$
%

\vspace{3ex}

\vspace{3ex}

\noindent {\bf Example 3.} A Caplet Price
\vspace{1ex}

The payoff at time $T'$ of a caplet for a LIBOR rate over the period $[T,T']$ can be written as [BGM, Section 3]
$$ C_{T'}=(T'-T)^{-1} (\frac{1}{P(T,T')} - \kappa)^+$$
Here  $\kappa = 1+ k (T'-T)$ and $k$ is the payment `cap'.  Since the monetary payment (paid at the end of the LIBOR period) is set at the beginning of the period, it is simple to compute the market value, $C_T$, of the caplet at time $T$.  The market value is determined by a fixed number of $T'$-maturity bonds.  An investment consisting of  $C_{T'} $ many discount bonds at time $T$ replicates the caplet payoff.  So, we see that
$$ C_{T}=P(T, T')(T'-T)^{-1}(\frac{1}{P(T,T')} - \kappa)^+$$
$$=(T'-T)^{-1}(1 - \kappa P(T, T'))^+$$

As in examples 1 and 2, the stopping times $\tau_n$ defined by equation (\ref{time}) allow each process
$$ \frac{P(t\wedge \tau_n,\tilde T)}{ P(t\wedge \tau_n,T_h)} $$
to be a martingale if $\tilde T = T'$ or $T$. The first process is dominated by the process with maturity parameter $T$ and both processes are bounded local martingales.

Following the pattern used to compute prices in previous examples, we consider a different contingent claim that expires at the smaller of the times $T$ and $\tau_n$. The payoff , $Cl_n$, for this claim is
$$Cl_n=(P(T\wedge \tau_n,T) - \kappa P(T\wedge \tau_n, T'))^+$$
For convenience we omit the constant factor of $(T' - T)^{-1}$ in defining this claim.  Clearly,  the limiting value of the payoff as $n\to \infty$ is the caplet value at time $T$ when multiplied by the constant factor.

We consider conditional expectations of the discounted payoff: 

\begin{equation}\label{capexp}
E[\  Cl_n /P(T\wedge \tau_n,T_h) \ | \ \mathcal{F}_t\ ]
\end{equation}
\vspace{2ex}
The event $D_n$ is defined by
\begin{equation*}
D_n = \{P(T\wedge \tau_n,T) > \kappa P(T\wedge \tau_n, T')\}
\end{equation*}
The expected value in equation (\ref{capexp}) can be written as

\begin{equation*}
E[ \ \frac{P(T\wedge \tau_n,T)}{P(T\wedge \tau_n,T_h)} 1_{D_n} \ | \ \mathcal{F}_t\ ] -\kappa E[\  \frac{P(T\wedge \tau_n,T')}{P(T\wedge \tau_n,T_h)} 1_{D_n} \ | \ \mathcal{F}_t\ ]
\end{equation*}

Since each integrand contains a discounted bond factor, we may use Bayes' Rule [KS,  Lemma 5.3]  to compute each conditional expectation with a change of measure. The process
$$\tilde \Lambda_t = \frac{P(t\wedge \tau_n, \tilde T)}{ P(t\wedge \tau_n,T_h)}\exp (-F(\tilde T)) $$
is a positive martingale whose expected value equals one for either choice $\tilde T = T'$ or $T$.  $\tilde \Lambda_T$ determines a probability measure $\tilde Q_n$ where
$$\frac{d\tilde Q_n}{dP} =\tilde  \Lambda_T$$ 
and Bayes' rule expresses the conditional expectation in equation (\ref{capexp}) as
\begin{equation}\label{measures}
\begin{aligned}
& \Lambda_t E_{Q_n}[\  1_{D_n}  \  |\  \mathcal{F}_t\ ] - \kappa   \Lambda_t' E_{ Q_n'}[\  1_{D_n}  \  |\  \mathcal{F}_t\ ]\\
&\ \   \\
&=\frac{P(t\wedge \tau_n,T)}{P(t\wedge \tau_n , T_h)}Q_n[\  D_n  \ |\ \mathcal{F}_t\ ]
-\kappa  \frac{P(t\wedge \tau_n,T')}{P(t\wedge \tau_n , T_h)} Q_n'[ \ D_n \  | \ \mathcal{F}_t\ ]\\
\end{aligned}
\end{equation}

In Appendix \ref{change} we show that

\begin{equation}\label{conprob1}
\begin{aligned}
 \lim_{n\to\infty} Q_n[\  D_n  \ |\ \mathcal{F}_t\ ]&= \mathcal{K}(T-t,\log \big[\frac{P(t,T)}{P(t , T')}\big], \log \kappa)\\
\mbox{where \hskip.1in} \mathcal{K}(s,x,y)&= Pr\{  \frac{2M_s}{2x^{-1} -A_s} \ge y\ \mbox{\ or\ }\ 2x^{-1} \le A_s\  \}\\
\end{aligned} 
\end{equation}
and
\begin{equation}\label{conprob2}
\begin{aligned}
 \lim_{n\to\infty} Q'_n[\  D_n  \ |\  \mathcal{F}_t\ ]&= \mathcal{K}'(T-t,\log \big[\frac{P(t,T)}{P(t , T')}\big], \log  \kappa)\\
\mbox{where \hskip.1in} \mathcal{K}'(s,x,y) & = Pr\{  \frac{2M_s}{2x^{-1} +A_s} \ge y\}\\
\end{aligned}
\end{equation}
Therefore, we can take limits in equation (\ref{measures}) to obtain the following price formula for the value of a caplet prior to the beginning of a LIBOR period:

\vspace{2ex}

\begin{equation}\label{caplet}
C_t=\delta^{-1}P(t,T)\mathcal{K}(\tau,\log \big[\frac{P(t,T)}{P(t , T')}\big], \log \kappa)
-\delta^{-1}\kappa P(t,T') \mathcal{K}'(\tau,\log \big[\frac{P(t,T)}{P(t , T')}\big], \log \kappa)
\end{equation}
\vspace{3ex}

\noindent where $\delta = T' - T$ and $\tau = T-t$.

\vspace{3ex}

\section{the volatility parameter}

\vspace{3ex}
 Each discounted bond price in Definition  \ref{dis} is a positive local martingale.  We will generalize this definition in a simple, but useful way.  For any $\sigma > 0$ it is obvious from Proposition \ref{arbitrage} that each process
\begin{equation}\label{speed}
\Big(\frac{P(0,T)}{P(0,T_h)}\Big)^{\big[\frac{2M_{\sigma^2 t}}{2+A_{\sigma^2 t}\log [P(0,T) /P(0,T_h)]}\big]}
\end{equation}
is a positive local martingale.  That is, if we adjust the filtration by defining
$$\mathcal{G}_t = \mathcal{F}_{\sigma^2 t}$$
then each process in equation (\ref{speed}) is a local martingale.  By speeding up (or slowing) the time evolution of asset prices with the factor $\sigma$ we introduce a volatility parameter within the model.  This is analogous to the appearance of a volatility parameter in the widely used log normal model for a stock price.  

In fact, there is no need to rescale the time in order to include this desirable feature.  The well known distributional equivalence of the processes $\{\sigma B_t\}$ and $\{B_{\sigma^2 t}\}$ allows us to modify  $M_t$ and $A_t$ so that processes are effectively speeded up or slowed down.

\vspace{3ex}

\begin{definition}\label{general} For $\sigma >0$ let $M_t^{(\sigma)} $ denote the simple exponential martingale
$$\exp (\sigma B_t - \sigma^2 t/2)$$
and let $A_t^{(\sigma)}$ denote its time integral
$$\int_0^t M_s^{(\sigma)}ds$$
For a given horizon date $T_h$ as in Definition \ref{horiz}, the $T_h$-maturity bond price of the \emph{general model} is given by 
\begin{equation}\label{generallong}
P(t,T_h) := \Big(\frac{P(0,t)}{P(0,T_h)}\Big)^{\big[\frac{-2M_t^{(\sigma)}}
{2+\sigma^2 A_t^{(\sigma)}\log [P(0,t) /P(0,T_h)]}\big]}
\end{equation}
For a maturity date $T < T_h$ the $T$-maturity bond price of  the \emph{general model} is given by 
\begin{equation}\label{general asset}  
P(t,T)= P(t,T_h)\Big(\frac{P(0,T)}{P(0,T_h)}\Big)^{\big[\frac{2M_t^{(\sigma)}}
{2+\sigma^2 A_t^{(\sigma)}\log [P(0,T) /P(0,T_h)]}\big]}
\end{equation}
for $t \leq T$.  The model with $\sigma = 1$ will be referred to as the \emph{special model}.
\end{definition}
\vspace{3ex}

\begin{proposition}  Each process
$$ \frac{P(t,T)}{P(t,T_h)} $$
in the general model, determined by the asset prices given in Definition \ref{general}, has the same distribution as the process in equation (\ref{speed}).  
\end{proposition}

\vspace{3ex}
\begin{proof}  A simple time change using Brownian scaling shows that the exponent process
$$  \frac{2M_t^{(\sigma)}}
{2+\sigma^2 A_t^{(\sigma)}\log (P(0,T) /P(0,T_h))}$$
of the general model discounted asset price has the same distribution (as a process) as the exponent 
$$\frac{2M_{\sigma^2 t}}{2+A_{\sigma^2 t}\log (P(0,T) /P(0,T_h))} $$
in equation (\ref{speed}).
\end{proof}
\vspace{3ex}

\begin{corollary}\label{gencap}  The caplet price formula in the general model for a payment cap of size $k$ on a LIBOR rate over the period $[T,T']$ is

\begin{equation}\label{generalcaplet}
\begin{aligned}
C_t=&\delta^{-1}P(t,T)\mathcal{K}(\sigma^2 \tau,\log \big[\frac{P(t,T)}{P(t , T')}\big], \log \kappa)\\
&-\delta^{-1}\kappa P(t,T') \mathcal{K}'(\sigma^2 \tau,\log \big[\frac{P(t,T)}{P(t , T')}\big], \log \kappa)\\
\end{aligned}
\end{equation}
\vspace{1ex}

\noindent In this equation $\delta = T' - T$, $1+\delta k = \kappa$, and $\tau = T-t$.

\end{corollary}

\vspace{3ex}
\begin{proof} The proposition implies that the special model asset price process, parametrized by $\sigma^2 t$, has the same distributional properties as those of the general model.  Then equation (\ref{caplet}), with the time rescaled, gives a replicating portfolio value for the general model.  
\end{proof}

\vspace{3ex}

\begin{corollary}  An approximate value of the  caplet price in the general model, valid  for small values of $\sigma$ or small values of $\tau = T-t$, the time until the LIBOR rate is fixed, is given by
\begin{equation}\label{approxcaplet}
\begin{aligned}
&C_t\approx \delta^{-1}P(t,T)\Phi (b_1) - \delta^{-1}\kappa P(t,T') \Phi (b_2)\\
\mbox{where \hskip.3in}\\
& b_i:=\frac{\log [P(t,T)/P(t , T')] -\log\kappa}{\sigma \sqrt{\tau}\log [P(t,T)/P(t , T')]}\ 
 \pm\  \frac{\sigma}{2}\sqrt{\tau}\log [P(t,T)/P(t , T')]
\\
\mbox{and \hskip.3in}&\phantom{b_2 = b_1 - \sigma \sqrt{\tau}\log [P(t,T)/P(t , T')]}\\
\end{aligned}
\end{equation}
$\Phi$ denotes the standard normal distribution function.
\end{corollary}

\vspace{3ex}

\begin{proof} The function $\mathcal{K}'(s,x,y)$ is defined in Example 3 of Section 4, and
$$\mathcal{K}'(\sigma^2\tau,x,y)= Pr\{  \frac{2M_{\sigma^2\tau}}{2x^{-1} +A_{\sigma^2\tau}} \ge y\}$$ 
which equals
$$Pr\{  \frac{2M_{\tau}^{(\sigma)}}{2x^{-1} +\sigma^2 A_{\tau}^{(\sigma)}} \ge y\}$$
The process $Y_{\tau}$, where
$$Y_{\tau} =  \frac{2M_{\tau}^{(\sigma)}}{2x^{-1} +\sigma^2 A_{\tau}^{(\sigma)}}\ ,$$
is determined by the SDE
$$dY_s = \sigma Y dB - \frac{\sigma^2}{2} Y^2 ds$$
and its initial condition $Y_0 = x$.  For small $\sigma$ we approximate $Y_t$ by $X_t$ where $X$ solves
$$dX_s = \sigma Y_0 dB - \frac{\sigma^2}{2} Y_0^2 ds$$
so that 
$$  X_s =x + \sigma xB_s - \frac{\sigma^2}{2}x^2s$$
Then
$$\mathcal{K}'(\sigma^2\tau ,x,y)\approx Pr\{x+\sigma x\sqrt{\tau}Z - \frac{\sigma^2}{2}x^2\tau \ge y\}$$
where $Z$ denotes a standard normal. Then
 \begin{equation*}
 \begin{aligned}
 \mathcal{K}' \approx & Pr\{x+\sigma x\sqrt{\tau}Z - \frac{\sigma^2}{2}x^2\tau \ge y\}\\
 =&Pr\{\sigma \sqrt{\tau}Z - \frac{\sigma^2}{2}x\tau \ge  \frac{y-x}{x}\}\\
 =&Pr\{Z \ge  \frac{y-x}{\sigma \sqrt{\tau}x} + \frac{\sigma}{2}\sqrt{\tau}x \}\\
 =& \Phi (\frac{x-y}{\sigma \sqrt{\tau}x} - \frac{\sigma}{2}\sqrt{\tau}x )\\
 \end{aligned}
 \end{equation*}

The other $\mathcal{K}$-function is approximately
$$\mathcal{K}(\sigma^2\tau,x,y) \approx Pr\{  \frac{2M_{\tau}^{(\sigma)}}{2x^{-1} -\sigma^2 A_{\tau}^{(\sigma)}} \ge y\}$$ 
For this case we use the SDE
$$dZ_s = \sigma x dB + \frac{\sigma^2}{2} x^2 ds$$
to approximate the process appearing in the definition of $\mathcal{K}$ . 
$$\mathcal{K} \approx  Pr\{x +\sigma x\sqrt{\tau}Z + \frac{\sigma^2}{2}x^2\tau \ge y\}$$

We substitute  the values 
$$x=  \log \big[\frac{P(t,T)}{P(t , T')}\big] \mbox{\  and\  }y = \log [\kappa] $$
to obtain approximate values for the functions  $\mathcal{K}$ and  $\mathcal{K'}$ in the caplet formula in Corollary (\ref{gencap}).  This gives an approximate caplet price. 
\end{proof}

\vspace{3ex}

\section{forward rates and the spot rate}\label{spotrate}

\vspace{3ex}

\begin{definition} A term structure admits a family of a forward rate processes, $r(t,T)$, if for each maturity date $T\le T_h$,
\begin{equation}\label{frates}
P(t,T) = \exp \big ( -\int_t^T r(t,u)du\big)
\end{equation}
for $t \le T$.
\end{definition}

\begin{proposition}\label{exist} The special model admits forward rate processes if and only if the nominal forward rate curve $F(t)$ satisfying  Definition \ref{curve} is an absolutely continuous function of   \ $t$ on the interval $[0, T_h]$.

\vspace{3ex}

Moreover, if $F(t)$ is absolutely continuous, so that 
$$ F(t) = \int_t^{T_h} f(s)ds$$
then for $T < T_h$
\begin{equation}\label{forward}
r(t,T) = \frac{f(T) }{M_t F(T)^2} \big[\log [P(t,T)/P(t,T_h)]\big]^2
\end{equation}
\end{proposition}

\vspace{3ex}

\begin{proof} Suppose that $r(t,T)$ is a family of forward rates. If $T < T_h$ we see  from equation (\ref{frates}) that
\begin{equation*}
P(t,T)/ P(t,T_h) = \exp \big ( \int_T^{T_h} r(t,u)du\big)
\end{equation*}
Then it follows from Proposition \ref{eee}  that
$$ \int_T^{T_h} r(t,u)du = \frac{2M_t }{2 F(T)^{-1} + A_t} $$
This equation may be solved for the function $F(T)$ and if the equation holds merely for some $t$ on an event with positive probability, the equality implies that $F(T)$ is the integral of an $L^1$ function.  Therefore, it is necessary that the nominal forward rate curve be absolutely continuous.

On the other hand, suppose that $F(t) = \int_t^{T_h} f(s)ds$ for some integrable function $f(s)$.  Then from Proposition \ref{eee}  we obtain
$$\log [P(t,T)/P(t,T_h)] = \frac{2M_t }{2 F(T)^{-1} + A_t}$$
$$ =\frac{2M_t }{2 [\int_T^{T_h} f(s)ds]^{-1} + A_t}$$
which implies that
$$\log P(t,T) =\frac{2M_t }{2 [\int_T^{T_h} f(s)ds]^{-1} + A_t}+\log P(t,T_h)$$
This identity shows that the process $\log P(t,T)$ is an absolutely continuous function of its maturity parameter and that it is the integral of
$$-\frac{4M_t }{(2 F(T)^{-1} + A_t)^2}\frac{f(T)}{F(T)^2}$$
That is,
$$r(t,T)=\frac{4M_t ^2}{(2 F(T)^{-1} + A_t)^2}\frac{f(T)}{M_tF(T)^2}$$
This establishes equation (\ref{forward}).  The forward rates also represent the long bond.  The formula above may be expressed  as
\begin{equation}\label{alternate}
r(t,T)=\frac{4f(T)M_t }{(2 + F(T)A_t)^2}
\end{equation}
Then, noting that $f = - F'$, we have
\begin{equation*}
\begin{aligned}
 \int_t^{T_h} r(t,u)du = &\frac{4M_t }{A_t(2 + F(u)A_t)}\Big ]_t^{T_h}\\
=&\frac{2M_t }{A_t} -  \frac{4M_t }{A_t(2 + F(t)A_t)}\\
=&\frac{2F(t)M_t }{2 + F(t)A_t} \\
\end{aligned}
\end{equation*}
and this is the formula for $-\log P(t,T_h)$ given in Proposition \ref{eee}.
\end{proof}

\vspace{3ex}

\begin{corollary}\label{P:forward} If the nominal forward curve is absolutely continuous, then for $t\le T< T_h$, each forward rate process $r(t,T)$ satisfies the SDE
\begin{equation}\label{E:forwardsde}
dr(t,T) = r(t,T) dB-r(t,T) \int^{T_h}_T r(t,u)du\cdot dt 
\end{equation}
\end{corollary}

\begin{proof} One may apply Ito's formula to the expression $r(t,T)$ in equation (\ref{alternate}) and verify that the SDE holds.
\end{proof}

\begin{corollary}\label{P:short}Suppose that the nominal forward curve is continuously differentiable on $[0,T_h)$. A continuous spot rate process $r(t)$ may be defined as 
\begin{equation}\label{spot}
r(t)=\frac{4f(t)M_t }{(2 + F(t)A_t)^2} = \frac{f(t) }{M_t F(t)^2} \big(\log [P(t,T_h)]\big)^2
\end{equation}
\end{corollary}

\begin{proof} In \cite{HJM} a spot rate $r(t)$ is defined within an HJM model as $\lim_{T\downarrow t}r(t,T)$.  If the function $f(t)$ is continuous, it follows from equation \ref{alternate} that
$$\lim_{T\downarrow t}r(t,T) = \lim_{T\downarrow t}\frac{4f(T)M_t }{(2 + F(T)A_t)^2}   =\frac{4f(t)M_t }{(2 + F(t)A_t)^2}$$ 
\end{proof}

\begin{remark}  The existence of a spot rate process is merely an artifact of the special model.  The model does not contain a money market asset; therefore the spot rate does not represent the rate of return on any investment.

An assumption that the forward rate curve is smooth on the closed interval $[0,T_h]$ allows us to define the \emph{longest forward rate} as the limit as $T\to T_h$ of $r(t,T)$.  It follows from equation (\ref{alternate}) that
$$ r(t,T_h) = f(T_h)M_t$$
Notice that the longest forward rate process is a geometric Brownian motion.
\end{remark}

\vspace{3ex}
\noindent{\bf Example.}   The \emph{Cantor function}, C(t), defined for $0 \le t \le 1$, is a well knonwn example of a non-absolutely continuous function [H].  $C(t)$ is monotone increasing on the interval $[0,1]$, $C(0) = 0$, and $C(1)=1$.  We let $T_h = 1$ and define a nominal forward curve function as
$$F(t) = 1 - C(t)$$
Then the special model determined by $F(t)$ does not admit forward rate processes.

In contrast, one choice for a smooth forward curve function is
$$F(t) = a(T_h - t)^b$$
where $a$ and $b$ are positive constants.  The function $f(t) = ab(T_h - t)^{b-1}$ is integrable even in the case that $b < 1$; for this choice of $b$ the initial instantaneous forward rate function, $f(t)$, is convex.  The special model determined by $F(t)$ has the forward rates processes (see equation (\ref{alternate}))
\begin{equation*}
r(t,T)=\frac{4abM_t }{(T_h - T)^{1-b}(2 + a(T_h - T)^{b}A_t)^2}
\end{equation*}
These rates blow up as the maturity parameter approaches the horizon date.  This is an example with forward rate processes, but no forward rates to the horizon.

\vspace{3ex}

\noindent{\bf Observations of instantaneous rate functions.\ }   The choice of a model horizon date may be  a matter of convenience.  On the other hand, one might set the model by selecting $T_h$ so that the corresponding long bond is actively traded and so that nearby maturities are also actively traded.  This choice for $T_h$ would allow a large data set to be collected in a brief time period. One might use this to set the \emph{expected value} of the longest forward rate
$$r(t,T_h).$$
It follows from the remark before the example that this expected value equals $f(T_h)$.  Therefore, one might choose a nominal rate function $F(t)$ which has a matching derivative at $t=T_h$,  $f(T_h)=-F'(T_h)$.

While observing the longest forward rate, one might also observe the spot rate.  The stochastic portion of the spot rate in equation (\ref{spot}) includes the same geometric Brownian motion that defines the stochastic behavior of the longest forward rate.  From (\ref{spot}) we obtain
\begin{equation*}
\begin{aligned}
r(t)= & \frac{f(t) }{M_t F(t)^2} \log ^2[P(t,T_h)]\\
= & \frac{f(T_h)f(t) }{r(t,T_h) F(t)^2} \log^2 [P(t,T_h)]\\
\end{aligned}
\end{equation*}
This implies that
\begin{equation}\label{nearzero}
\frac{r(t)r(t,T_h)}{\log ^2[P(t,T_h)]} =  \frac{f(T_h)f(t) }{F(t)^2} 
\end{equation}

The identity in equation (\ref{nearzero}) states that the deterministic function 
$$f(T_h)f(t)/F(t)^2$$
is determined by values of the long bond, the longest forward rate, and the spot rate.  This is one example where a modeler might choose the function $f(t)$, at least for $t$ near zero, to approximate observations of $r(t)$, $r(t,T_h)$, and $P(t,T_h)$.
\vspace{3ex}

\section{origin of the model}

A risk-neutral measure in an HJM model [HJM, section 5] , [MR, equation (13.20)] produces the following dynamics for forward rate functions $r(t,T)$:
$$ dr = \sigma(t,T) \int_t^T \sigma(t,u)du \cdot dt + \sigma (t,T) dB$$
The choice of 
$$  \sigma (t,T) = r(t,T)$$
is considered in [HJM , equation (45)].  The dynamics become
$$ dr = r\cdot  \int_t^T r(t,u)du \cdot dt + r dB$$
Morton \cite{Mo} found that this equation allows interest rate explosions. Consequently, this model allows arbitrage.   Explosions occur because of  the quadratic behavior of the drift term.  A positive initial value, combined with a volatility which prevents the rate from becoming negative, allows uncontrolled growth to infinity.

Another dynamics for forward interest rates within an HJM model is given in [MR2, example 2.1]:
$$ dr = -\sigma(t,T) \int_T^{T_h} \sigma(t,u)du\cdot dt + \sigma (t,T) dB$$
In this situation, when we specialize to $  \sigma (t,T) = r(t,T)$ we obtain
\vspace{3ex}
$$ dr = -r(t,T) \int_T^{T_h} r(t,u)du\cdot  dt + r(t,T) dB$$
\vspace{3ex}

\noindent The appearance of the negative sign in the (quadratic) drift is due to a change of measure.  This is the dynamics when a forward measure corresponding to the time $T_h$ is chosen.  With such a measure,  at least formally, each asset price discounted by the long bond price is a martingale.

Solutions of this equation do not explode.  And, conveniently, the solutions are explicit.  If we state this equation in terms of the quantity
$$ Y_t =  \int_T^{T_h} r(t,u)du$$
the equations reads
$$ dY = - \frac{1}{2}Y^2 dt + Y dB$$
This is the logarithm of the process we define in Definition \ref{dis}.

\vspace{3ex}

\section*{Acknowledgements}

\vspace{1ex}

The authors wish to thank Michael Tehranchi, Liheng Xu, and Siyu Zhang for helpful conversations.  We also thank a referee for suggestions concerning an earlier version of this paper.

\vspace{3ex}


\appendix

\section{ }\label{change} We identify two conditional probabilities that appear in Example 3, where the price formula for a caplet is obtained.   In Example 3 we used stopping times  $\tau_n$ defined by equation (\ref{time}) where $0 < T <T' \le T_h$. For $t \le T$, the process
$$\tilde \Lambda_t = \frac{P(t\wedge \tau_n, \tilde T)}{ P(t\wedge \tau_n,T_h)}\exp (-F(\tilde T)) $$
is a positive martingale whose expected value equals one for either choice $\tilde T = T'$ or $T$.  $\tilde \Lambda_T$ determines a probability measure $\tilde Q_n$ where
$$\frac{d\tilde Q_n}{dP} =\tilde  \Lambda_T$$ 
For each $n = 1, 2, \dots $ we let 
\begin{equation*}
D_n = \{\frac{P(T\wedge \tau_n,T)}{P(T\wedge \tau_n, T')} > \kappa \}
\end{equation*}
where $\kappa$ is a positive constant.

\vspace{2ex}
\begin{proposition}

\begin{equation}\label{conprobA}
\begin{aligned}
\lim_{n\to\infty} Q_n[\  D_n  \ |\ \mathcal{F}_t\ ]&= 
Pr\{  \frac{2M_{T-t}}{2x_t^{-1} -A_{T-t}} \ge \log \kappa \ \mbox{or} \   A_{T-t} \ge 2x_t^{-1}\}\\
\mbox{where \hskip.8in} x_t&=\frac{P(t,T)}{P(t, T')} \\
\mbox{and \hskip1in}&   \\
\lim_{n\to\infty} Q'_n[\  D_n  \ |\ \mathcal{F}_t\ ]&= 
Pr\{  \frac{2M_{T-t}}{2x_t^{-1} +A_{T-t}} \ge \log \kappa \}\\
\end{aligned}
\end{equation}

\end{proposition}

\vspace{3ex}
\begin{remark}  The probability terms in the statement of the Proposition involve an independent process $(M_s, A_s)$.  That is, distribution functions of the random variables 
$$  \frac{2M_{s}}{2x^{-1} \pm A_{s}}$$
are used to express the limit.
\end{remark}

\vspace{3ex}

\begin{proof}  We first consider the case $\tilde T = T'$.  Then
$$ \Lambda_T' = \frac{P(T\wedge \tau_n, T')}{ P(T\wedge \tau_n,T_h)}\exp (-F(T')) $$
With respect to the measure $Q_n'$, the process $P(s,T)/P(s, T')$ is a Markov process.  In order to see this, we use Proposition \ref{eee} to find that
$$\log \big(  \frac{P(s, T)}{ P(s,T')}\big)= \frac{2M_{s}}{2F(T)^{-1} +A_{s}} -
\frac{2M_{s}}{2F(T')^{-1} +A_{s}}$$
which we write as
\begin{equation}\label{log}
\log \big(  \frac{P(s, T)}{ P(s,T')}\big)= Y_s - Y'_s
\end{equation}
Moreover, from Proposition \ref{eee} we also have
$$ \Lambda_T' = \exp (Y'_{T\wedge \tau_n}-F(T')) $$
Now both processes $Y_s$ and $Y'_s$ satisfy the SDE in equation (\ref{sde}).  In particular, 
$$ dY_s' = Y_s'dB_s - \frac{1}{2}(Y_s')^2 ds$$
This provides a very useful formula for $\Lambda'$:
$$ \Lambda_T' = \exp (\int_0^{T\wedge \tau_n} Y'_sdB_s -\frac{1}{2}\int_0^{T\wedge \tau_n} (Y'_s)^2 ds ) $$
This shows that $\Lambda '$ is a stochastic exponential determined by the bounded drift process 
$Y'_s1_{s\le \tau_n}$.  We apply the Girsanov theorem to see that $\Lambda'$ adds the drift $Y'_s$ to the Brownian motion.  That is, with respect to the measure $Q'_n$ the process
\begin{equation}\label{drift'}
B'_s = B_s - \int_0^{s\wedge \tau_n} Y'_udu
\end{equation}
is a standard Brownian motion.  Now we derive an SDE for the process in equation (\ref{log}).  First, we have
\begin{equation}\label{difference}
d(Y - Y') = (Y - Y')dB - \frac{1}{2}( Y^2 - (Y')^2)ds
\end{equation}
since each process satisfies equation (\ref{sde}).  However, under $Q'_n$ the process $B_s$ is a semimartingale and we express the stochastic integral using equation (\ref{drift'}):
$$d(Y - Y') = (Y - Y')(dB'+Y'ds) - \frac{1}{2}( Y^2 - (Y')^2)ds$$
which holds path-wise for $s \le \tau_n$.  We obtain
\begin{equation}\label{difference'}
\begin{aligned}
d(Y - Y') = &(Y - Y')dB' - \frac{1}{2}(-2YY'+2(Y')^2+ Y^2 - (Y')^2)ds\\
= &(Y - Y')dB' - \frac{1}{2}(Y - Y')^2ds\\
\end{aligned}
\end{equation}
The autonomous form of this SDE shows that $(Y_s - Y'_s)$ is a Markov process.

One immediate consequence is that
$$Q'_n[\  D_n  \ |\ \mathcal{F}_t\ ] = Q'_n[\  D_n  \ |\ \log [\frac{P(t\wedge \tau_n,T)}{P(t\wedge \tau_n, T')}] \ ]$$
However, a further consequence of equation (\ref{difference'}) is that we can obtain an explicit description of $Y_s - Y_s'$:
\begin{equation}\label{desc}
Y_s - Y_s' = \frac{2 M'_{s}}{2(F(T)-F(T'))^{-1} + A'_{s}}
\end{equation}
In this equation, $M'$ denotes the elementary exponential martingale for the Brownian motion $B'$ and $A'$ denotes its time integral.  Now the event $D_n$ is
$$D_n = \{ Y_{T\wedge \tau_n} - Y'_{T\wedge \tau_n} \ge \log \kappa\}$$
and so we have
\begin{equation*}
Q'_n[\  D_n  \ |\ \mathcal{F}_t\ ] = Q'_n[\  Y_{T\wedge \tau_n} - Y'_{T\wedge \tau_n} \ge \log \kappa  \ |\ Y_{t\wedge \tau_n} - Y'_{t\wedge \tau_n}\ ] 
\end{equation*}
Since $\tau_n$ is a stopping time for the process in equation (\ref{desc}), we use the transition probabilities for the process $Y-Y'$ to compute the conditional probability. Denoting
$$Pr\{    \frac{2 M'_{s}}{2x^{-1} + A'_{s}}\in E\}$$
by $p(s,x,E)$, we apply [KS, Proposition 6.6, Chapter 2]:

\begin{equation}\label{transition}
Q'_n[\  D_n  \ |\ \mathcal{F}_t\ ] = p(T\wedge \tau_n - t\wedge \tau_n,Y_{t\wedge \tau_n} - Y'_{T\wedge \tau_n}, [\log \kappa, \infty)] 
\end{equation}
Equation (\ref{transition}) expresses the conditional probability in terms of a fixed function whose only dependence on $n$ is through the stopping time $\tau_n$.

Since $\tau_n \to \infty$ almost surely as $n\to\infty$, we have
\begin{equation*}
\lim_{n\to\infty}Q'_n[\  D_n  \ |\ \mathcal{F}_t\ ] = p(T-t,Y_t-Y'_t, [\log \kappa, \infty))
\end{equation*}
This establishes the second limit in the Proposition.

\vspace{3ex}

We next compute $Q_n[\  D_n  \ |\ \mathcal{F}_t\ ] $ using a similar method.  We have 
$$ \Lambda_T = \exp (Y_{T\wedge \tau_n}-F(T)) $$
and since the process $Y$ satisfies the SDE given in equation (\ref{sde}),
$$ \Lambda_T = \exp (\int_0^{T\wedge \tau_n} Y_sdB_s -\frac{1}{2}\int_0^{T\wedge \tau_n} Y_s^2 ds ) $$
so that the measure $Q_n$ adds the drift $Y_s1_{s\le \tau_n}$ to the Brownian motion.  With respect to $Q_n$ the process
\begin{equation*}
\tilde B_s = B_s - \int_0^{s\wedge \tau_n} Y_u du
\end{equation*}
is a standard Brownian motion.  From equation (\ref{difference}) we obtain
\begin{equation*}
\begin{aligned}
d(Y - Y') = &(Y - Y')(d\tilde B+Yds) - \frac{1}{2}( Y^2 - (Y')^2)ds\\
= &(Y - Y')d\tilde B + \frac{1}{2}( Y-Y')^2ds\\
\end{aligned}
\end{equation*} 
Again, the SDE for the difference process is autonomous, but now there is an explosion.  The SDE has a positive quadratic drift term. An explicit solution is
\begin{equation}\label{explosion}
Y_s - Y_s' = \frac{2 \tilde M_{s}}{2(F(T)-F(T'))^{-1} - \tilde A_{s}}
\end{equation}
We denote the transition probabilities for this Markov process,
$$Pr\{    \frac{2 \tilde M_{s}}{2x^{-1} - \tilde A_{s}}\in E\ \mbox{\ or\ }  \tilde A_s \ge 2x^{-1}\},$$
by $q(s,x,E\cup \{\infty\})$. Then Proposition 6.6 of [KS] applies and we find that
$$Q_n[\  D_n  \ |\ \mathcal{F}_t\ ] = q(T\wedge \tau_n - t\wedge \tau_n,Y_{t\wedge \tau_n}-Y'_{t\wedge \tau_n}, [\log \kappa, \infty])$$ 
As before, the only dependence on $n$ is through the stopping times and we obtain
\begin{equation*}
\lim_{n\to\infty}Q_n[\  D_n  \ |\ \mathcal{F}_t\ ] = q(T-t,Y_t-Y'_t, [\log \kappa, \infty])
\end{equation*}

\end{proof}

\vspace{3ex}

\begin{remark} The exploding property of the process $Y_s - Y'_s$ in equation (\ref{explosion}) reveals that changing the numeraire can result in a singular measure.  The measure $Q_n$ is defined in terms of the positive martingale
$$\Lambda (t) =\frac{P(t\wedge \tau_n,T) }{P(t\wedge \tau_n,T_h)}\cdot \frac{P(0,T_h) }{P(0,T)} $$ 
Suppose that the equation 
$$\frac{d Q}{dP} :=\frac{1 }{P(T,T_h)}\cdot \frac{P(0,T_h) }{P(0,T)}$$ 
did define a probability measure $Q$.  Such a measure would correspond, formally, to a change of numeraire where discounted asset prices of the form
$$\frac{P(t,S) }{P(t,T)}$$
would be local martingales.  However, for maturity dates $S$ exceeding $T$,  the price $P(t,S)$ would  be zero with positive $Q$-measure.  Equation (\ref{explosion}) would give an explicit formula for
$$Y'_s - Y_s = \log \frac{P(t,T') }{P(t,T)}$$
assuming that 
\begin{equation*}
\tilde B_s = B_s - \int_0^{s} Y_u du
\end{equation*}
were a standard Brownian motion.  But then the event  $Y'_s - Y_s =-\infty$ would correspond to $(F(T) - F(T'))\tilde A_s\ge 2$.  Since $Y'_s - Y_s$ is finite with probability one in the model, the measure $Q$ could not be absolutely continuous with respect to the measure $P$.
\end{remark}
\vspace{3ex}


\section{ }\label{pitfall}

\vspace{3ex}

A naive use of martingale theory to compute prices for contingent claims may lead to incorrect answers.  Here is one example where a forward contract price for the long bond can not be computed in the customary manner. The exercise price required to receive a long bond is denoted by $\kappa$ and the expiration date is  $T < T_h$.  

One might believe that the discounted contract price is given by
\begin{equation}\label{longcon}
E[ (P(T,T_h) - \kappa)/ P(T,T_h)\ |\  \mathcal{F}_{t} \ ]
\end{equation}
$$  = 1 - \kappa E[  \frac{P(T,T)}{P(T,T_h)} \ |\  \mathcal{F}_{t} \ ] $$
Theorem \ref{positive} asserts that the process $P(t,T) / P(t,T_h)$ is a positive local martingale, and it is well known that such processes are also supermartingales.  Therefore,
$$E[  \frac{P(T,T)}{P(T,T_h)} \ |\  \mathcal{F}_{t} \ ] \leq P(t,T) / P(t,T_h)$$
In this case, the inequality is strict.  To see this, suppose that equality holds with probability one.  Then we have
$$E[ \  \frac{P(T,T)}{P(T,T_h)} \ ] = E[ \ P(t,T) / P(t,T_h)\ ]$$
Proposition \ref{eee} allows us to write this identity as
$$E[ \  \exp (\frac{2M_T }{2 F(T)^{-1} + A_T}) \ ] = E[ \ \exp (\frac{2M_t }{2 F(T)^{-1} + A_t})\ ]$$
However, in [ GK, Corollary 4.2] we show that the expected value on the left hand side is strictly less than the other expected value.  By using  the Markov property of the process
$$Y_t = \frac{2M_t }{2 F(T)^{-1} + A_t}$$
one can prove the stronger result that  
$$ E[  P(T,T) / P(T,T_h) \ |\  \mathcal{F}_{t} \ ] < P(t,T) / P(t,T_h)$$
with probability one.  Therefore, the computed price from  equation (\ref{longcon}) is larger than the correct value
$$P(t,T_h) - \kappa P(t,T)$$
\

\vspace{3ex}

\vspace{3ex}


\begin{thebibliography}{DGMS}

\bibitem[BGM]{BGM} A. Brace, D. Gatarek, and M. Musiela (1997): \emph{The market model of interest rate dynamics}, Math. Fin. \textbf{7}, 127-147.

\bibitem[Do]{Do} M. Dothan  (1990): \emph{Prices in Financial Markets}, Oxford University Press.

\bibitem[Du]{Du} D. Dufresne (2001): \emph{The integral of geometric Brownian motion}, Adv. in Appl. Probab. \textbf{33}, 223-241.

\bibitem[DH]{DH} P. Dybvig and C. Huang (1989): \emph{Nonnegative wealth, absence of arbitrage, and feasible consumption plans}, Review Finan. Studies \textbf{1}, 377-401.

\bibitem[GK]{GK} V. Goodman and K. Kim (preprint):  \emph{Exponential martingales and the time integral of geometric Brownian motion}.

\bibitem[KS]{KS} I. Karatzas and S.Shreve (1991): \emph{Brownian Motion
and Stochastic Calculus}, Springer-Verlag New York.

\bibitem[WH]{WH} B. Wong and C.C. Heyde (2004): \emph{On the martingale property of stochastic exponentials}, J. Appl. Probab. \textbf{41}, 654-664.

\bibitem[H]{H} P. R. Halmos (1974): \emph{Measure Theory}, Springer-Verlag New York.

\bibitem[HJM] {HJM} D. Heath, R. Jarrow, and A. Morton (1992): \emph{Bond
pricing and the term structure of interest rates}: A new methodology for contingent claim valuation, Econometrica \textbf{60}, 77-105.

\bibitem[Mi]{Mi} K. Miltersen (1994): \emph{An arbitrage theory of the
term structure of interest rates}, Ann. Appl. Prob. \textbf{4},
953-967

\bibitem[MSS]{MSS} K. Miltersen, K. Sandmann, and D. Sondermann (1997): \emph{Closed form solutions for term structure derivatives with log-nornal interest rates}, J. Finance \textbf{LII},
409-430

\bibitem[Mo]{Mo} A.J. Morton (1989): \emph{Arbitrage and martingles}
Doctoral dissertation, Cornell University, Ithaca, U.S.A.

\bibitem[MR] {MR} M. Musiela and M. Rutkowski (1997): \emph{Martingale
Methods in Financial Modelling}, Applications of Mathematics
\textbf{36}, Springer-Verlag, NewYork

\bibitem[MR2]{MR2} M. Musiela and M. Rutkowski (1997):
\emph{Continuous-time term structure models: Forward measure
approach}, Finance stochast. \textbf{1}, 261-291

\bibitem[P] {P} P. Protter (1990): \emph{Stochastic Integration and Differential Equations}, Applications of Mathematics \textbf{21}, Springer-Verlag, NewYork

\end{thebibliography}
\end{document}